# Solving the two envelopes problem with the Intermediate Amount Strategy

Tsikogiannopoulos Panagiotis

**Abstract**

This paper introduces a strategy in the two envelopes problem that utilizes the prior beliefs of two players about the amount of money that their envelopes can contain. This strategy gives them more information about the decision of switching they have to make when one of the envelopes is opened and the amount it contains is revealed. The player who implements this strategy can predict which amount is larger with a probability greater than 1/2 so that he will have a greater expected return from the game than that of the other player who will not use the same strategy. The symmetrical case in which both players implement the same strategy is also analyzed.

**Introduction**

The strategy that is presented below is applied to the two envelopes problem played by two players with the following rules: Each player is given an envelope containing an amount of money. Both players are informed that one of the envelopes contains double the amount of the other envelope without knowing which envelope contains which amount. The envelope of a player is opened and the amount it contains is revealed. Then each player, secretly and independently from the other one, has to decide whether to ask or not for the exchange of his envelope with that of the other player. If both players request for an exchange then it will be made, otherwise they will keep their original envelopes. Finally, each player will win the amount of money inside the envelope that ended up to him.

The Intermediate Amount Strategy (IAS) can be applied only in the variation of the game where one of the two amounts is revealed and only if the player who will implement it has a prior belief about the two amounts the game is played with. The cases where the players have no prior beliefs are well analyzed in the existing literature [1], [2], [3].

Let's see how this strategy works: First, the player thinks of an amount $M$. The optimal method of determining $M$ will be discussed below, but the strategy works with any amount $M$. Then the amount contained in one envelope is revealed to the player. Assume that this amount is $P$. If $P > M$, the player will assume that $P$ is the larger of the two amounts contained in the envelopes. If $P \leq M$, the player will assume that $P$ is the smaller of the two amounts. So according to his assumption, he will decide whether to keep or ask for the exchange of his envelope in order to end up with the larger amount.

Applying this strategy leads to the following conclusions: In the case where $M$ is smaller than the smaller amount of the two envelopes then the probability of the player to ask for the envelope containing the larger amount is 1/2. In the case where $M$ is larger than the larger amount of the two envelopes then the probability of the player to ask for the envelope containing the larger amount is also 1/2. However, in the case where $M$ lies between the amounts of the two envelopes then the probability

of the player to ask for the larger amount is 1. It becomes clear that, however small is
the probability of the amount *M* to be set between the amounts of the two envelopes,
the IAS gives a probability of a correct prediction always greater than 1/2.
Furthermore, the more realistically a player sets *M*, the greater his probability
becomes to choose the larger amount.
If both players implement equally successfully the IAS then the strategy's advantage
is neutralized, but each player does not know whether the other one will apply the IAS
and how successfully he will apply it. As we will see, it turns out that it is profitable
for each player to apply the IAS to maximize his expected amount, regardless of what
the other player does.
A general analysis of the strategy is described in the paper of Dov Samet, Iddo Samet
and David Schmeidler [4].
Let's see how we can formalize and quantify the profit obtained using the IAS. We
will examine two ways that a player can implement the IAS, according to the
estimations he can make about the two amounts.

## 1. 1st Type of calculation – Uniform Distribution

In the 1st type of calculation, the player considers the distribution of the probability of
the two amounts lying within a particular interval to be uniform. More specifically,
the player will choose this type of calculation if both of the following conditions are
true:
A) He has a prior belief for the maximum amount of money that the organizer of the
    game could pay to a player. So he sets in a variable *N* the maximum amount that
    according to his estimation an envelope can contain.
B) He cannot estimate where the two amounts are placed within the interval $(0, N]$
    that he set.

Definitions

At first we will specify the rules by which the game is played. We consider that the
game organizer has decided beforehand the two amounts that he will place in the
envelopes. He places the smaller amount in one envelope and the larger amount in the
other one and with a random process he selects which envelope to give to each player.
We will refer to the unknown amounts of the two envelopes with the continuous
variables *X* and *2X* that will correspond to the smaller and the larger amount
respectively. The domains of the variables *X* and *2X* according to the player who set *N*
are the following:

$$X \in (0, N/2] \qquad\qquad 2X \in (0, N] \qquad\qquad (1.1)$$

The reason for defining the maximum *X* as half of *N* is that if the player believed that
*X* could be even larger, then the resulting *2X* would have been greater than *N*, so he
would have to adjust the value of *N*.
In the 1st type of calculation the player considers the probability densities of the
variables *X* and *2X* to be constant, i.e. it is equally likely for *X* and *2X* to lie anywhere
within their domains and that they cannot be found elsewhere. We define the constant
probability density function *f*(*s*) to be the inverse of the product of the total range of
*N* and the percentage *s* of this range that we will refer to. Formally:

$$f(s) = \frac{1}{sN} \tag{1.2}$$

By this definition, the probabilities that the two variables $X$ and $2X$ have to lie within their domains are respectively:

$$P(0 < X \le N/2) = \int_0^{N/2} f(1/2)dx = 1 \qquad P(0 < 2X \le N) = \int_0^N f(1)dx = 1 \tag{1.3}$$

The mean amounts of the two envelopes are respectively:

$$\overline{X} = \int_0^{N/2} x \cdot f(1/2)dx = \frac{N}{4} \qquad \overline{2X} = \int_0^N x \cdot f(1)dx = \frac{N}{2} \tag{1.4}$$

The formulas (1.3) and (1.4) come from the definition of the Uniform Distribution. Now, for the definition of the intermediate amount $M$, we saw in the Introduction that the player has an interest to choose $M$ to the point that maximizes the probability to lie between the amounts of $X$ and $2X$. We can also say that a given $M$ lies between $X$ and $2X$ when $X$ takes a value between $M/2$ and $M$. So the probability of $M$ to lie between the two amounts is maximized when the range between $M/2$ and $M$ is also maximized and this happens when the player sets $M$ to the maximum value that $X$ can get, namely when:

$$M = \frac{N}{2} \tag{1.5}$$

In this case, $M$ will lie between the two amounts when $X$ lies between $N/4$ and $N/2$ and the probability of this to happen is given by the quotient of the favorable range over the total range that $X$ can be found at. Formally:

$$P(X \le M \le 2X) = \frac{(N/2) - (N/4)}{N/2} = \frac{1}{2} \tag{1.6}$$

So, the player who will apply the 1st type of calculation should always set the value $N/2$ to $M$ to consider that there is a probability of $1/2$ for $M$ to lie between the amounts of the two envelopes.
We will denote by $A$ the amount that both players see in the envelope of player A. We will denote by $M_A$ the intermediate amount $M$ that is set by player A.

Calculations

The calculations of the players must be made at the time where one of the envelopes is opened in order to have the maximum possible information about the game.
We initially consider that each player doesn't know what strategy the other player is about to implement.
Before the envelope is opened, the expected amount contained in player A's envelope is given by assigning a probability of $1/2$ to the events of holding the smaller or the larger amount, namely:

$$E_{init}[A] = \frac{1}{2}X + \frac{1}{2} \cdot 2X = \frac{3X}{2} \tag{1.7}$$

The following calculations concern player A when his envelope is opened.
We distinguish two cases depending on the relationship between the amount $A$ that is revealed and the amount $M_A$ he set:

1. If $A \leq M_A$ then according to the IAS player A will consider that the amount of $A$ contained in his envelope is the smaller of the two amounts and therefore he will request to exchange his envelope.

   We now distinguish two subcases depending on whether $A$ is actually the smaller or the larger amount:

   1a. If $A = X$ then provided that the exchange request will be accepted, player A will end up with the amount $2X$.
   1b. If $A = 2X$ then provided that the exchange request will be accepted, player A will end up with the amount $X$.

According to the formulas (1.1) and (1.5) player A considers that the subcase 1a is twice more likely to occur than the subcase 1b. This is because $P(2X \leq N/2) = P(X \leq N/4) = \frac{1}{2} P(X \leq N/2)$ due to the uniform distribution. Since subcases 1a and 1b are complementary, player A has to assign them the probabilities of 2/3 and 1/3 respectively.
So the expected amount given to player A in case 1, provided that his exchange request will be accepted, is:

$$E_{accepted}[A] = \frac{2}{3} \cdot 2X + \frac{1}{3}X = \frac{5X}{3} \tag{1.8}$$

If player A's exchange request is not accepted then his expected amount is calculated by alternating $X$ with $2X$ in formula (1.8), i.e.:

$$E_{denied}[A] = \frac{2}{3}X + \frac{1}{3} \cdot 2X = \frac{4X}{3} \tag{1.9}$$

2. If $A > M_A$ then according to the IAS player A will consider that the amount of $A$ contained in his envelope is the larger of the two amounts and therefore will ask to keep his envelope.

In case 2, player A considers the event that his envelope containing the amount $X$ to be impossible, according to formulas (1.1) and (1.5).
So he considers certain that the expected amount in his envelope is:

$$E_{keep}[A] = 2X \tag{1.10}$$

We should note that in case 2 the expected amount of player A is not depended on player B's decision to exchange or not his envelope because player A decides to keep his envelope and that cannot be changed.

We will now calculate the expected amount given by the IAS for player A before his envelope is opened, i.e. by co-calculating cases 1 and 2:
The probability of $2X$ to be laid in the interval $(0, N/2]$ is the same with the probability to be laid in the interval $(N/2, N]$ due to the uniform distribution.
We showed above that the probability of $X$ to be laid in the interval $(0, N/2]$ is double the probability of $2X$ to be laid in the same interval. Hence, the probability of case 1 to occur is 3 times larger than the probability of case 2 to occur and since these two events are complimentary, the corresponding probability values are 3/4 and 1/4 respectively.

We will denote by $p_B$ the probability that player A assigns to the event that player B will ask to exchange his envelope. Then the expected amount that player A calculates it will end up to him, before his envelope is opened, is given by the formula:

$$E[A] = \frac{3}{4}\left[p_B \cdot E_{accepted}[A] + (1 - p_B) \cdot E_{denied}[A]\right] + \frac{1}{4} E_{keep}[A] = \left(\frac{3}{2} + \frac{p_B}{4}\right) X \qquad (1.11)$$

Even if player A is certain that player B will not ask for an exchange ($p_B = 0$), the result of the formula (1.11) becomes the same as that of the formula (1.7), which is his expected amount without the use of IAS. The larger player A considers $p_B$ to be, the larger his expected amount becomes. In the ideal case where player B knows nothing about IAS, he might always agree for the exchange of envelopes because he would believe that it doesn't make a difference anyway. If player A knows this in advance, he can set $p_B = 1$ in formula (1.11) and the result will be the maximum possible expected amount that is equal to $7X / 4$. If player A has no information about player B's intention to exchange or not his envelope then he will not be able to set a value to probability $p_B$, but even in this case, formula (1.11) shows that it is in player A's best interest to implement the IAS.

Of course, player A may be wrong about his estimations. For example it may be the case that the actual amount $X$ is larger than the amount $M_A$ he set. In this and every case that $M_A$ is not set between the amounts $X$ and $2X$, player A will not have a benefit by implementing the IAS and his expected amount will be $3X / 2$ according to formula (1.7). Nevertheless, player A considers that there is a probability of 1/2 for $M_A$ to actually be set between $X$ and $2X$ according to the formula (1.6) and then he will end up with the amount $2X$ if he will ask to keep his envelope or if his exchange request will be accepted. So, he can also calculate his expected amount by the formula:

$$E[A] = \frac{1}{2} \cdot \frac{3X}{2} + \frac{1}{2} \cdot 2X = \frac{7X}{4} \qquad (1.12)$$

which again gives the maximum expected amount for player A, given that if he asks for an exchange, his request will be accepted.

Suppose now that player B is also implementing the IAS with the 1$^{st}$ type of calculation and he considers his own uniform distribution. Therefore, he will set his own maximum amount that an envelope can contain, let's name it $N_B$, and he will set his own intermediate amount $M_B = N_B / 2$. The formulas (1.8), (1.9) and (1.10) remain the same for player B and they give him the expected amounts of player A with respect to $M_B$.

If player B wants to calculate his expected amount before player A's envelope is opened, he will accordingly set a probability $p_A$ that player A will ask for an exchange and his expected amount will be:

$$E[B] = \frac{3}{4} \cdot \frac{5X}{3} + \frac{1}{4}[p_A \cdot 2X + (1-p_A) \cdot X] = \left(\frac{3}{2} + \frac{p_A}{4}\right)X \quad (1.13)$$

Note that the results of the formulas (1.11) and (1.13) differ only to the probability of the other player to ask for an exchange, as it was expected due to the symmetry of the two cases.

We will now examine the expected amounts when both players implement the IAS, but with different intermediate amount $M$ each. We will rely on the results of the formulas (1.8), (1.9) and (1.10) and we will examine every possible interrelation between the amounts $A$, $M_A$ and $M_B$ and every possible outcome of the exchange request. The results are shown in the following table:

**Expected amounts when both players implement IAS**

| | Interrelation between $A, M_A, M_B$ | E[A] when player A's exchange request: | | | E[B] when player B's exchange request: | | |
|---|---|---|---|---|---|---|---|
| | | is accepted | is not accepted | player A isn't requesting | is accepted | is not accepted | player B isn't requesting |
| a. | $M_A < M_B < A$ | | | 2X | | X | |
| b. | $M_B < M_A < A$ | | | 2X | | X | |
| c. | $A < M_A < M_B$ | | $\frac{4X}{3}$ | | | | $\frac{5X}{3}$ |
| d. | $A < M_B < M_A$ | | $\frac{4X}{3}$ | | | | $\frac{5X}{3}$ |
| e. | $M_A < A < M_B$ | | | 2X | | | $\frac{5X}{3}$ |
| f. | $M_B < A < M_A$ | $\frac{5X}{3}$ | | | 2X | | |

*Table 1*

In the cases a, b, c, d of Table 1, an exchange is not agreed because when one player is requesting for an exchange according to the IAS the other player also according to the IAS chooses to keep his envelope. We also note that the total expected amount of the players in these cases is 3X, because the amounts they are going to share are the X and the 2X.

In the case e, none of the players according to IAS will ask for an exchange because they both consider it more likely to have the amount of 2X in their envelope. This disagreement is due to the fact that each M is positioned to a different side in relation to the amount A.

In the case f, both players according to the IAS will ask for an exchange and it will be granted, because they both consider it more likely to have the amount of *X* in their envelope. This disagreement is also due to the fact that each *M* is positioned to a different side in relation to *A*, but reversed from the case e.
In the cases e and f, the total expected amount of the players results greater than 3*X* and this is because only the player who has set *M* at the correct side of *A* is right.

Now, let's calculate the expected amount of each player, before the amount *A* is revealed, when they both know beforehand that the other player will also implement the IAS. We will use the results of Table 1 excluding the cases e and f, because in these cases one of the two players has made a wrong assessment. We have already mentioned that the cases a and b have a probability to occur equal to 1/4 and the cases c and d have a probability to occur equal to 3/4. Thus, the expected amounts derived from Table 1 are:

$$E_{both\ IAS}[A] = \frac{1}{4} \cdot 2X + \frac{3}{4} \cdot \frac{4X}{3} = \frac{3X}{2} \tag{1.14}$$

$$E_{bothIAS}[B] = \frac{1}{4} X + \frac{3}{4} \cdot \frac{5X}{3} = \frac{3X}{2} \tag{1.15}$$

The results of the formulas (1.14) and (1.15) are the same, due to the symmetry of the cases. It is obvious that when both players know in advance that the other player will also implement the IAS they can expect no profit from the strategy, unless they are in the cases e and f. Of course they will also have to implement the IAS because otherwise they will give the advantage to the other player who will.

## 2. 2<sup>nd</sup> Type of calculation – Normal Distribution

In the 2<sup>nd</sup> type of calculation, the player considers that the probability of the two amounts to deviate from the values he set is following a normal distribution. More specifically, the player will choose this type of calculation if both of the following conditions are true:
A) He doesn't want to estimate the maximum amount of money that the organizer of the game could pay to a player.
B) He has a prior belief of how large the two amounts can be.

Definitions

In the 2<sup>nd</sup> type of calculation the variable *N* is not defined.
Player A will denote by $X_A$ what he considers to be the most likely smaller amount and thus he will denote by $2X_A$ what he considers to be the most likely larger amount. These two definitions will be made before the opening of an envelope.
Player A will also set the amount $M_A$ in the middle between $X_A$ and $2X_A$, i.e.:

$$M_A = \frac{3X_A}{2} \tag{2.1}$$

The amount $M_A$ is considered by player A to have the greatest probability to lie between the actual amounts of the two envelopes, which are the $X$ and the $2X$, similarly to the way is set in formula (1.5) of the 1st type of calculation.

Here, the probability density function $f$ of $X_A$ to be any amount is not constant as of the 1st type of calculation, but it has a peak to the point that the player considers as the most likely smaller amount and decreases as we deviate from it. Under this consideration, the most proper way to define $f$ is the truncated normal distribution function. So let's define its parameters.

We define the standard deviation $\sigma$ as a function of a point $\mu$ where $f$ has its peak, i.e.:

$$\sigma(\mu) = \mu \cdot c_v \tag{2.2}$$

The reason we set the standard deviation as a function depended by the variable $\mu$ is to assure that the certainty of the player that has set correctly an amount is depended by the value of that amount. In other words, it is not equally likely for a player to deviate by 100 euros when he predicts that the smaller amount is 0.01 euros and when he predicts that it is 1000 euros. The deviation in the first case should be smaller and in the second case should be larger.

$c_v$ is the coefficient of variation and it remains constant. Its value will be determined below.

Each time we are referring to the probability densities of $X_A$ and $2X_A$, we will set $\mu = X_A$ and $\mu = 2X_A$ accordingly.

We are now ready to define the probability density function, following the theory of the truncated normal distribution [5], [6], when the lower bound of the independent variable is zero and the upper bound is infinity:

$$\varphi(t) = \frac{1}{\sqrt{2\pi}} e^{-\frac{t^2}{2}} \tag{2.3}$$

$$\Phi(x) = \int_{-\infty}^{x} \varphi(t) dt \tag{2.4}$$

$$f(x, \mu) = \frac{\frac{1}{\sigma(\mu)} \varphi\left(\frac{x - \mu}{\sigma(\mu)}\right)}{1 - \Phi\left(\frac{-\mu}{\sigma(\mu)}\right)} \tag{2.5}$$

The functions $\varphi(t)$ and $\Phi(x)$ are respectively the density and the cumulative distribution of the standard normal distribution N(0,1). The $f(x, \mu)$ is the probability density function that we will use in our calculations. Its arguments are the variable $x$ that corresponds to the amount for which the probability is calculated and the constant $\mu$ that corresponds to the point that the function has its peak. We set the additional argument of $\mu$ to the function $f$ so that we can use the same function for both the amounts $X_A$ and $2X_A$.

To calculate the coefficient of variation $c_v$ that we left undefined above, we note that exactly at the middle between the amounts $X_A$ and $2X_A$, the probability densities of the two amounts must be equal. Formally:

$$f(M_A, X_A) = f(M_A, 2X_A) \tag{2.6}$$

This equation should be true because for example if the player considers that the most likely smaller amount is 100 euros, so the most likely larger amount is 200 euros, then given that the envelope is opened contains 150 euros, he will have to assign equal probabilities to the events that the 150 euros is the smaller or the larger amount. Solving the equation (2.6) by the variable $c_v$ we get:

$$c_v = \frac{3}{4\sqrt{6}\ln 2} \approx 0.3678 \tag{2.7}$$

The domains of $X_A$ and $2X_A$ are:

$$X_A \in (0, \infty) \qquad\qquad 2X_A \in (0, \infty) \tag{2.8}$$

The probabilities of the two amounts to lie within their domains are respectively:

$$P[0 < X_A < \infty] = \int_0^\infty f(x, X_A)dx = 1 \qquad P[0 < 2X_A < \infty] = \int_0^\infty f(x, 2X_A)dx = 1 \tag{2.9}$$

The value of $f(2X_A, 2X_A)$ is half the value of $f(X_A, X_A)$ and in general the graph of $f(x, 2X_A)$ is shorter and wider than that of $f(x, X_A)$. This is due to the fact that when the player is off by say 10% in the value of $X_A$ over the actual $X$ then he is necessarily off by 10% in the value of $2X_A$ over the actual $2X$. The same must be true with the corresponding probabilities of $X_A$ and $2X_A$ to lie within two specific intervals that are set relatively to these values. Formally, the function $f$ has been defined so that the following equation is always true for any percentage $k$:

$$\int_{X_A - k \cdot X_A}^{X_A + k \cdot X_A} f(x, X_A)dx = \int_{2X_A - k \cdot 2X_A}^{2X_A + k \cdot 2X_A} f(x, 2X_A)dx \tag{2.10}$$

The above equation is ensured by the definition of the standard deviation as a function of point $\mu$ in the formula (2.2) and it is independent to $c_v$.

Again, we will denote by $A$ the amount that is revealed in the envelope of player A. In Graph 1 is presented an example of plotting the functions $f(x, X_A)$ and $f(x, 2X_A)$ for $X_A = 300$, $2X_A = 600$, $M_A = 450$ and a revealed amount $A = 550$.

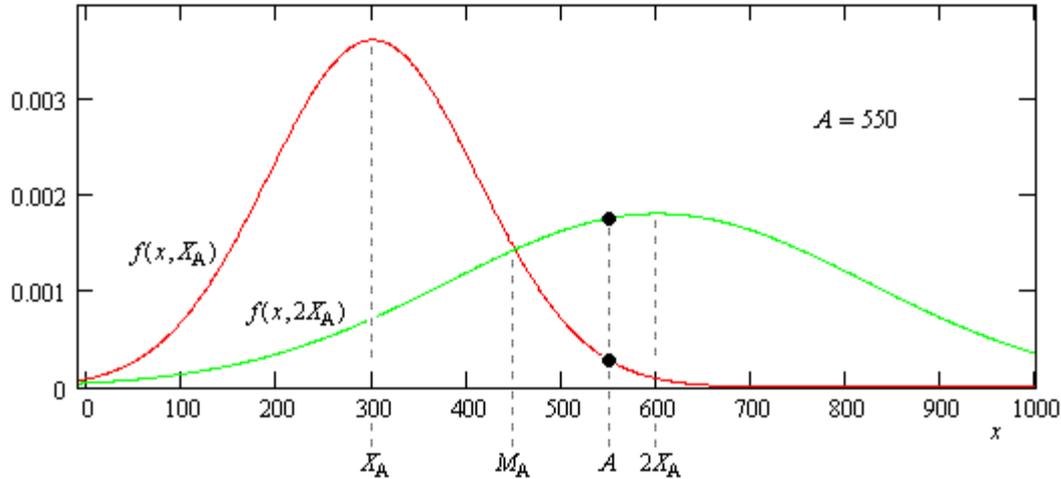

*Graph 1*

Calculations

The following calculations are concerning player A when his envelope is opened. When player A's envelope is opened he will implement the IAS. Even though it is equally likely for player A to hold the smaller or the larger amount, when his amount is revealed he can utilize his prior belief for the values of the two amounts, so that he can attribute a different probability to the events that the amount he sees is the smaller or the larger one.

First, player A will calculate the density functions $f(A, X_A)$ and $f(A, 2X_A)$ and then he will calculate the final probabilities of the amount $A$ he sees to be the smaller or the larger one, in the following way:

$$P[A = X] = \frac{f(A, X_A)}{f(A, X_A) + f(A, 2X_A)} \tag{2.11}$$

$$P[A = 2X] = \frac{f(A, 2X_A)}{f(A, X_A) + f(A, 2X_A)} \tag{2.12}$$

so that:

$$P[A = X] + P[A = 2X] = 1 \tag{2.13}$$

We now distinguish two cases depending on the relationship between the amount $A$ that is revealed and the amount $M_A$ that he set:

1. If $A \leq M_A$ then according to the IAS player A will ask for the exchange of his envelope. The expected amount that will end up to him, provided that his exchange request will be accepted, is:

$$E_{accepted}[A] = P[A = X] \cdot 2X + P[A = 2X] \cdot X \tag{2.14}$$

If his exchange request will not be accepted then his expected amount is calculated by switching $X$ with $2X$ in formula (2.14):

$$E_{denied}[A] = P[A = X] \cdot X + P[A = 2X] \cdot 2X \tag{2.15}$$

2. If $A > M_A$ then according to the IAS player A will ask to keep his envelope, so the expected amount that will end up to him is:

$$E_{keep}[A] = P[A = X] \cdot X + P[A = 2X] \cdot 2X \tag{2.16}$$

By the way we defined the probability density function, or simply by observing Graph 1, it follows that:

$$\forall A < M_A \rightarrow P[A = X] > P[A = 2X] \tag{2.17}$$

$$\forall A > M_A \rightarrow P[A = 2X] > P[A = X] \tag{2.18}$$

By the formulas (2.14), (2.16), (2.17) and (2.18) it follows that for every initial amount $A$, player A considers that the probability he has to end up with the amount $2X$ is always greater than the probability to end up with the amount $X$. So the implementation of the IAS is also advantageous for the player in the 2$^{nd}$ type of calculation.

If player B wants to implement the IAS, he will have to set as $X_B$ the amount that he considers as the most likely to be the smaller amount and as $2X_B$ the amount that he considers as the most likely to be the larger one. He will also have to set his own amount $M_B$ as the middle between $X_B$ and $2X_B$ and his own function $f$ for the probability densities of the two amounts.

In a similar way as that we showed in the 1$^{st}$ type of calculation, the two players can calculate their expected amounts before the opening of the envelope and when they both know that the other player will implement the IAS also.

**Conclusions**

We analyzed two types of implementing the Intermediate Amount Strategy (IAS) that utilize the prior beliefs of two players about the amounts that their envelopes can contain. These methods suggested whether or not the players should ask for the exchange of their envelopes. The decision of the players for which type of calculation should implement must not depend on which type gives them the larger expected amount, but on which of the two distributions' parameters think they can define with greater accuracy. The more accurate are defined the parameters of each type of calculation the larger the expected amount will be, but the opposite is not true. We saw that in both the 1$^{st}$ and the 2$^{nd}$ type of calculation the player has an interest to implement the IAS because he raises his chances to end up with the larger amount, even with the least accurate definition of its parameters as opposed to the indifferent way of choosing an envelope. The Intermediate Amount Strategy gives a new interest to the old two envelopes problem and shows how many things can still be written about it.